\documentclass[MSNbibl,number,citesort,seceqn,dvips]{arxbj}
\usepackage{mathrsfs}

\aid{0}
\volume{19}
\issue{5A}
\pubyear{2013}
\firstpage{1776}
\lastpage{1789}
\doi{10.3150/12-BEJ429} 

\makeatletter

\newproclaim{definition}{Definition}[section]
\newtheorem{lemma}[definition]{Lemma}
\newtheorem{theorem}[definition]{Theorem}
\newremark{remark}[definition]{Remark}
\newtheorem{corollary}[definition]{Corollary}
\newremark{example}[definition]{Example}
\newcommand{\eqref}[1]{(\ref{#1})}
\def\P{ {\mathsf P}}
\def\E{ {\mathsf E}}
\def\oo{\infty}
\makeatother

\begin{document}
\begin{frontmatter}

\title{Some inequalities of linear combinations of independent random variables: II}
\runtitle{Linear combinations of independent random variables}

\begin{aug}
\author[1]{\fnms{Xiaoqing} \snm{Pan}\thanksref{1,e1}\ead[label=e1,mark]{panxq@mail.ustc.edu.cn}},
\author[2]{\fnms{Maochao} \snm{Xu}\thanksref{2}\ead[label=e2]{mxu2@ilstu.edu}} \and
\author[1]{\fnms{Taizhong} \snm{Hu}\corref{}\thanksref{1,e3}\ead[label=e3,mark]{thu@ustc.edu.cn}}
\runauthor{X. Pan, M. Xu and T. Hu} 
\address[1]{Department of Statistics and Finance, School of Management, University of Science and Technology of China, Hefei, Anhui
230026, People's Republic of China.\\ \printead{e1}; \printead*{e3}}
\address[2]{Department of Mathematics, Illinois State University, Normal, IL 61761, USA.\\ \printead{e2}}
\end{aug}

\received{\smonth{2} \syear{2011}}
\revised{\smonth{1} \syear{2012}}

\begin{abstract}
Linear combinations of independent random variables have been extensively studied in the
literature. However, most of the work is based on some specific distribution assumptions. In this
paper, a companion of (\textit{J.~ Appl. Probab.} \textbf{48} (2011) 1179--1188), we unify the study of linear combinations of independent
nonnegative random variables under the general setup by using some monotone transforms. The
results are further generalized to the case of independent but not necessarily identically
distributed nonnegative random variables. The main results complement and generalize the results
in the literature including (In \textit{Studies in Econometrics, Time Series, and Multivariate Statistics} (1983) 465--489 Academic Press;
 \textit{Sankhy\={a} Ser. A} \textbf{60} (1998) 171--175; \textit{Sankhy\={a} Ser. A} \textbf{63} (2001) 128--132; \textit{J. Statist. Plann. Inference} \textbf{92} (2001) 1--5;
  \textit{Bernoulli} \textbf{17} (2011) 1044--1053).
\end{abstract}


\end{frontmatter}

\section{Introduction}

Linear combinations of independent nonnegative random variables arise naturally in statistics,
operations research, reliability theory, computer science, economic theory, actuarial science and
other fields. There are a large number of extensive studies on this topic in the literature. Some
typical applications could be found in \cite{1,2,10,15,22} and references therein. It should be remarked that most
of the work in the literature is under some specific distribution assumptions such as exponential,
Weibull, Pareto and gamma, etc.

Under the general framework, Karlin and Rinott \cite{8} studied the linear combinations of nonnegative
independent and identically distributed (i.i.d.) random variables $X_1,X_2,\ldots, X_n$ with $X^p_i$ having a
log-concave density for $0<p<1$. They showed that if $q<0$ and $p^{-1}+q^{-1}=1$, then
\begin{equation}\label{krst}
    (a^q_1,\ldots, a^q_n)\,{\succeq}_\mathrm{m}\, (b^q_1,\ldots, b^q_n)
      \quad \Longrightarrow  \quad\sum_{i=1}^n a_i X_i\ge_\mathrm{st}  \sum_{i=1}^n b_i X_i,
\end{equation}
where $\mathbf{a}=(a_1,\ldots, a_n)\in \Re^n_+$, $\mathbf{b}=(b_1,\ldots, b_n)\in \Re^n_+$, $\succeq_\mathrm{m}$
means the majorization order, and $\ge_\mathrm{st}$ means the usual stochastic
order (their formal definitions are given in Section~\ref{s2}).

Recently, Yu \cite{22} further studied this problem and showed two very interesting results. For
nonnegative i.i.d. random variables $X_1,X_2,\ldots,X_n$, if $\log(X_i)$ has a log-concave density, then
\begin{equation}\label{yu1}
    (\log a_1,\ldots, \log a_n) \succeq_\mathrm{m} (\log b_1, \ldots, \log b_n)
      \quad \Longrightarrow  \quad\sum_{i=1}^n a_i X_i\ge_\mathrm{st}  \sum_{i=1}^n b_i X_i.
\end{equation}
In contrast to the result \eqref{krst} in \cite{8}, Yu \cite{22} proved that if
$X^p_i$ has a log-concave density for $p>1$, then, for $q>1$ and $p^{-1}+q^{-1}=1$,
\begin{equation}\label{yu2}
  (a^q_1,\ldots, a^q_n) \, {\succeq}_\mathrm{m}  \, (b^q_1,\ldots, b^q_n)
     \quad \Longrightarrow \quad \sum_{i=1}^n a_i X_i\le_\mathrm{st}  \sum_{i=1}^n b_i X_i.
\end{equation}

However, in practical situation, it is quite often that random variables may not be i.i.d., that is,
i.i.d. seems to be a restrictive assumption; see \cite{9,11,12,23} and references therein. Xu and Hu \cite{21} successfully extended
\eqref{yu1} to the case of independent but not necessarily identically distributed random
variables under some suitable conditions.

This paper is a companion of \cite{21}. In this paper, we will further study this topic.
First, in Section \ref{s3}, we unify the study of \eqref{krst}--\eqref{yu2} by using some monotone
transforms, and then extend the results to the case of independent but not necessarily identically
distributed nonnegative random variables in Section \ref{s4}. Some examples are highlighted as well.

\section{Preliminaries}\label{s2}

In this section, we recall the definitions of some stochastic orders and majorization orders, which will be
used in the sequel.

\begin{definition}
Let $X$ and $Y$ be two random variables with distribution functions $F$ and $G$, density functions
$f$ and $g$ (if exist), respectively. Then $X$ is said to be smaller than $Y$

\begin{itemize}
\item in the usual stochastic order, denoted by $X \le_\mathrm{st} Y$, if $F(x)\ge G(x)$ for all
$x$;

\item in the likelihood ratio order, denoted by $X \le_\mathrm{lr} Y$, if
$g(x)/f(x)$ is increasing in $x$ for which the ratio is well defined.
\end{itemize}
\end{definition}

The likelihood ratio order is stronger than the usual stochastic order. For more discussions on
stochastic orders, please refer to \cite{18}.\eject

We shall also be using the concept of majorization in our discussion. For extensive and
comprehensive details on the theory of majorization orders and their applications, please refer to
\cite{14}. Let $x_{(1)}\le x_{(2)}\le\cdots \le x_{(n)}$ be the increasing
arrangement of components of the vector $\mathbf{x}=(x_1,x_2,\ldots,x_n)$.

\begin{definition}
For vectors $\mathbf{x}, \mathbf{y} \in \Re^n$,  $\mathbf{x}$ is said to be
\begin{itemize}
\item majorized by $\mathbf{y}$, denoted by $\mathbf{x}\preceq_\mathrm{m} \mathbf{y}$, if
$\sum_{i=1}^n x_{(i)}=\sum_{i=1}^n y_{(i)}$ and
\[
     \sum_{i=1}^j x_{(i)}\ge \sum_{i=1}^j y_{(i)}\qquad \mbox{for }
     j=1,\ldots,n-1;
\]

\item weakly supmajorized by  $\mathbf{y}$, denoted by $\mathbf{x}\preceq^\mathrm{w}\mathbf{y}$, if
\[
    \sum_{i=1}^j x_{(i)}\ge \sum_{i=1}^j y_{(i)} \qquad\mbox{for } j=1,\ldots,
    n;
\]

\item weakly submajorized by  $\mathbf{y}$, denoted by $\mathbf{x}\preceq_\mathrm{w}\mathbf{y}$, if
\[
    \sum_{i=j}^n x_{(i)}\le \sum_{i=j}^n y_{(i)} \qquad \mbox{for } j=1,\ldots, n.
\]
\end{itemize}
\end{definition}

A real-valued function $h$ defined on a set $A\subseteq \Re^n$ is said to be Schur-concave
[Schur-convex] on $A$ if, for any $\mathbf{x}, \mathbf{y}\in A$,
\[
   \mathbf{x}\succeq_\mathrm{m} \mathbf{y} \quad\Longrightarrow \quad h(\mathbf{x}) \le [\ge]\, h(\mathbf{y}),
\]
and $h$ is log-concave on $A$ if $A$ is a convex set and, for any $\mathbf{x}, \mathbf{y}\in A$ and
$\alpha\in [0,1]$,
\[
   h\bigl(\alpha \mathbf{x} + (1-\alpha) \mathbf{y}\bigr) \ge [h(\mathbf{x})]^\alpha [h(\mathbf{y})]^{1-\alpha}.
\]

To prove the main results in the next section, we recall the following two well-known lemmas. The
first one gives the preservation properties of the weakly majorization orders under monotone
transforms, while the second one states that the log-concavity is closed under integral.

\begin{lemma}[(\cite{14}, Theorem 5.\textrm{A}.2)]\label{le-20110722}
\begin{enumerate}[(iii)]
\item[(i)] For all increasing and convex functions $g$,
  \[
     \mathbf{x}\preceq_\mathrm{ w} \mathbf{y} \quad \Longrightarrow \quad(g(x_1), \ldots, g(x_n))\preceq_\mathrm{ w}
         (g(y_1), \ldots, g(y_n)).
  \]
\item[(ii)] For all increasing and concave functions $g$,
  \[
     \mathbf{x}\preceq^\mathrm{ w} \mathbf{y}  \quad\Longrightarrow \quad(g(x_1), \ldots, g(x_n))\preceq^\mathrm{ w}
         (g(y_1), \ldots, g(y_n)).
  \]
\item[(iii)] For all decreasing and convex functions $g$,
  \[
     \mathbf{x}\preceq^\mathrm{ w} \mathbf{y}  \quad\Longrightarrow \quad(g(x_1), \ldots, g(x_n))\preceq_\mathrm{ w}
         (g(y_1), \ldots, g(y_n)).
  \]
\item[(iv)] For all decreasing and concave functions $g$,
  \[
     \mathbf{x}\preceq_\mathrm{ w} \mathbf{y}  \quad\Longrightarrow \quad(g(x_1), \ldots, g(x_n))\preceq^\mathrm{ w}
         (g(y_1), \ldots, g(y_n)).
  \]
\end{enumerate}
\end{lemma}

\begin{lemma}[(\cite{4,17})]\label{le-logconcave} Suppose that $h\dvtx\Re^m\times \Re^k
\to \Re_+$ is a log-concave function and that
 \[
        g(\mathbf{x}) =\int_{\Re^k} h(\mathbf{x}, \mathbf{z})\, \mathrm{d} \mathbf{z}
 \]
 is finite for each $\mathbf{x}\in \Re^m$. Then $g$ is log-concave on $\Re^m$.
\end{lemma}

\section{i.i.d. nonnegative random variables}\label{s3}

In this section, we unify the study of linear combinations of i.i.d. nonnegative random variables
under the general setup by using some monotone transforms.

\begin{theorem}\label{th-20110209}
Let $X_1,X_2, \ldots, X_n$ be i.i.d. absolutely continuous and nonnegative random variables, and let
$\phi,\psi\dvtx \Re_+\to \Re_+$ be two twice continuously differentiable and strictly monotone
functions such that, for all $(u,v)\in\Re_+^2$,
\begin{equation}\label{eq-20110210-1}
      \phi''(u)\ge 0
\end{equation}
and
\begin{equation}\label{eq-20110210-2}
    \phi''(u)\psi''(v)\phi(u)\psi(v)\ge [\phi'(u)\psi'(v)]^2.
\end{equation}
Assume that $\psi^{-1}(X_1)$ has a log-concave density function, where $\psi^{-1}$ is the inverse
function of~$\psi$. If
\begin{equation}\label{eq-20110211-3}
    (\phi^{-1}(a_1),\ldots, \phi^{-1}(a_n) )\succeq_\mathrm{ m}  (\phi^{-1}(b_1),\ldots,
      \phi^{-1}(b_n) ),
\end{equation}
then
\begin{equation}\label{eq-20110721-3}
            \sum_{i=1}^n a_i X_i\ge_\mathrm{ st} \sum_{i=1}^n b_i X_i.
\end{equation}
Moreover, if $\phi$ is increasing (decreasing), the majorization order in \eqref{eq-20110211-3}
can be replaced by the submajorization (supmajorization) order.
\end{theorem}

\begin{pf} First, we prove that \eqref{eq-20110211-3} implies \eqref{eq-20110721-3}. To see it,
suppose that \eqref{eq-20110211-3} holds. Fix any $t\in \Re_+$, and define
\[
   h(\mathbf{c}) = \P\Biggl(\sum_{i=1}^n \phi(c_i) X_i\le t \Biggr).
\]
It suffices to show that $h(\mathbf{c})$ is Schur-concave in $\mathbf{c}\in \Re^n_+$. Define
\[
    A=\Biggl\{(\mathbf{y},\mathbf{c})\in \Re^{2n}_+ \dvtx \sum_{i=1}^n \phi(c_i) \psi(y_i)\le t\Biggr\}.
\]
Then
\[
  h(\mathbf{c}) = \P \Biggl(\sum_{i=1}^n \phi(c_i) \psi (\psi^{-1}(X_i) ) \le t \Biggr)
     = \int_{\Re^n} g(\mathbf{y}, \mathbf{c}) \, \mathrm{d} \mathbf{y},
\]
where
\[
   g(\mathbf{y}, \mathbf{c}) =1_{A}(\mathbf{y}, \mathbf{c})\cdot \prod_{i=1}^n f_{\psi^{-1}(X_i)}(y_i).
\]
Next, we will discuss when
\[
    \Phi(u,v)=\phi(u)\psi(v)
\]
is convex on $\Re^2_+$. Note that the Hessian matrix for $\Phi(u,v)$ is
\[
   \pmatrix{   \phi''(u)\psi(v) & \phi'(u)\psi'(v) \cr
           \phi'(u)\psi'(v) & \phi(u)\psi''(v)      }.
\]
It is known that if the Hessian matrix for $\Phi(u,v)$ is nonnegative semi-definite, then $\Phi(u,v)$ is
convex on $\Re^2_+$. That is, if \eqref{eq-20110210-1} and \eqref{eq-20110210-2} hold, then $\Phi(u,v)$ is
convex and, hence, $A$~is a convex set. This implies that $1_A(\mathbf{y},\mathbf{c})$ is log-concave on $(\mathbf{
y},\mathbf{c})\in\Re^{2n}$. Thus, $g(\mathbf{y},\mathbf{c})$ is log-concave. By Lemma \ref{le-logconcave}, $h(\mathbf{c})$ is log-concave. Since $h(\mathbf{c})$ is permutation symmetric, and the permutation symmetry and
log-concavity imply Schur-concavity (see Fact 3.1 in \cite{20}), we conclude that $h(\mathbf{c})$ is
Schur-concave.

Next, suppose that $\phi$ is decreasing and $ (\phi^{-1}(a_1),\ldots, \phi^{-1}(a_n)
)\succeq^\mathrm{ w}  (\phi^{-1}(b_1),\ldots, \phi^{-1}(b_n) )$. By Proposition 5.A.9 in
\cite{14}, there exists a real vector $(c_1,\ldots,c_n)\in\Re_+^n$ such that
\[
    (\phi^{-1}(a_1),\ldots, \phi^{-1}(a_n) ) \le (c_1, \ldots, c_n)
       \succeq_\mathrm{ m}  (\phi^{-1}(b_1),\ldots, \phi^{-1}(b_n) ).
\]
Here, for two vectors $\mathbf{s}, \mathbf{t}\in\Re^n$, $\mathbf{s}\ge \mathbf{t}$ means componentwise
ordering. Since $\phi$ is decreasing, we have $a_i\ge \phi(c_i)$ for each $i$ and, hence, $
    \sum^n_{i=1} a_i X_i \ge_\mathrm{ st} \sum^n_{i=1} \phi(c_i) X_i.
$
On the other hand, it is shown that
$
  \sum^n_{i=1} \phi(c_i) X_i \ge_\mathrm{ st} \sum^n_{i=1} b_i X_i.
$
Thus, we conclude \eqref{eq-20110721-3}.

Finally, suppose that $\phi$ is increasing and $ (\phi^{-1}(a_1),\ldots, \phi^{-1}(a_n)
)\succeq_\mathrm{ w}  (\phi^{-1}(b_1),\ldots,\break \phi^{-1}(b_n) )$. Again, by Proposition
5.A.9 in \cite{14}, there exists a real vector $(c_1,\ldots,c_n)\in\Re_+^n$ such
that
\[
    (\phi^{-1}(a_1),\ldots, \phi^{-1}(a_n) ) \ge (c_1, \ldots, c_n)
       \succeq_\mathrm{ m}  (\phi^{-1}(b_1),\ldots, \phi^{-1}(b_n) ).
\]
The rest of the proof is similar and is hence omitted. This completes the proof.
\end{pf}

\begin{theorem}\label{th-20110210}
Let $X_1,X_2, \ldots, X_n$ be i.i.d. absolutely continuous and nonnegative random variables, and let
$\phi, \psi\dvtx \Re_+\to \Re_+$ be two twice continuously differentiable and strictly monotone
functions such that, for all $(u,v)\in\Re_+^2$,
\begin{equation} \label{eq-20110210-3}
      \phi''(u)  \le 0
\end{equation}
and \eqref{eq-20110210-2} hold. If $\psi^{-1}(X_1)$ has a log-concave density function, then
\begin{equation}\label{eq-20110211-4}
    (\phi^{-1}(a_1),\ldots, \phi^{-1}(a_n) )\succeq_\mathrm{ m}  (\phi^{-1}(b_1),\ldots,
      \phi^{-1}(b_n) ) \quad\Longrightarrow \quad\sum_{i=1}^n a_i X_i\le_\mathrm{ st} \sum_{i=1}^n b_i X_i.
\end{equation}
Moreover, if $\phi$ is increasing (decreasing), the majorization order in \eqref{eq-20110211-4}
can be replaced by the supmajorization (submajorization) order.
\end{theorem}

\begin{pf} The proof is similar to that of Theorem \ref{th-20110209} by observing that, for any
$t\in \Re_+$,
\[
   h(\mathbf{c}) = \P\Biggl(\sum_{i=1}^n \phi(c_i) X_i> t \Biggr)
\]
is Schur-concave in $\mathbf{c}\in \Re^n_+$ under conditions \eqref{eq-20110210-2} and
\eqref{eq-20110210-3}.
\end{pf}

\begin{remark}
The results in Theorems \ref{th-20110209} and \ref{th-20110210} can be extended to permutation invariant
random variables; see \cite{13}. This was also pointed out by one of the referees.
\end{remark}

\begin{remark}
Theorems \ref{th-20110209} and \ref{th-20110210} do not apply to the case that $\phi(x)=x$ or/and
$\psi(x)=x$. One may wonder whether $\sum^n_{i=1} a_i X_i$ and $\sum^n_{i=1} b_i X_i$ are ordered
in the usual stochastic order whenever $\mathbf{a}, \mathbf{b}\in \Re_+^n$ such that $\mathbf{a}\succeq_\mathrm{ m} \mathbf{b}$ under the assumption that $X_i$ has a log-concave density. However,
this is not true. A counterexample is given by Diaconis and Perlman \cite{3} as follows: For $X_i$
having a gamma distribution with shape parameter $\alpha\ge 1$ (whose density is log-concave) and
$n\ge 3$, if $\mathbf{a}$ and $\mathbf{b}$ differ in exactly two components, then the distribution
functions of $\sum^n_{i=1} a_i X_i$ and $\sum^n_{i=1} b_i X_i$ are of unique crossing. In fact, if
$\mathbf{a}\succeq_\mathrm{ m} \mathbf{b}$, then $\E [\sum^n_{i=1} a_i X_i ] = \E
[\sum^n_{i=1} b_i X_i ]$ and, hence, there cannot be a stochastic order between the two
weighted sums unless there is equality in distribution.
\end{remark}

\begin{remark}\label{re-20110213}
For Theorem \ref{th-20110209}, conditions \eqref{eq-20110210-1} and \eqref{eq-20110210-2} imply
that $\phi$ and $\psi$ are both convex; while, for Theorem \ref{th-20110210},
\eqref{eq-20110210-2} and \eqref{eq-20110210-3} imply that $\phi$ and $\psi$ are both concave.
Some special choices of $\phi$ and $\psi$ in Theorem \ref{th-20110209} or \ref{th-20110210} are as
follows.

\begin{itemize}
\item $\phi(x)=\psi(x)=\mathrm{e^\mathit{x}}$:\\
Conditions \eqref{eq-20110210-1} and \eqref{eq-20110210-2} are satisfied. Theorem \ref{th-20110209}
reduces to Theorem 1 in \cite{22}. That is, if $\log X_i$ has a log-concave density, we have
\begin{equation}\label{eq-20110213-1}
   (\log a_1, \ldots, \log a_n) \succeq_\mathrm{ w} (\log b_1, \ldots, \log b_n)\quad\Longrightarrow
       \quad\sum^n_{i=1} a_i X_i \ge_\mathrm{ st} \sum^n_{i=1} b_i X_i.
\end{equation}

\item $\phi(x)=x^{1/q}$ and $\psi(x)=x^{1/p}$:\\
It can be checked that
\begin{eqnarray*}
    \eqref{eq-20110210-1} \mbox{ and } \eqref{eq-20110210-2}
    \mbox{ hold}
          \quad  & \Longleftrightarrow &\quad (p,q)\in A_0\cup A_1 \cup A_2;\\
  \eqref{eq-20110210-2} \mbox{ and } \eqref{eq-20110210-3} \mbox{ hold}\quad & \Longleftrightarrow &\quad (p,q)\in A_3,
\end{eqnarray*}
where
\begin{eqnarray*}
     A_0&=&  \{(p,q)\dvtx p<0,q<0 \},    \\
     A_1&=& \biggl \{(p,q)\dvtx p<0,0<q<1,\frac 1p+ \frac 1q\ge 1\biggr \},\\
    A_2&=& \biggl\{(p,q)\dvtx 0<p<1,q<0,\frac 1p+ \frac 1q\ge 1 \biggr\},\\
    A_3 &=& \biggl\{(p,q)\dvtx p>1,q>1,\frac 1p+ \frac 1q\le 1 \biggr\}.
\end{eqnarray*}
Choosing $(p,q)\in A_1, A_2$ and $A_3$, Theorems \ref{th-20110209} and \ref{th-20110210} reduce to
Corollaries \ref{co-20110212-1}, \ref{th-Karlin} and~\ref{th-Yu}, respectively. For more details,
see Remark \ref{re-20110721}. Choosing $(p,q)\in A_0$, from Theorem~\ref{th-20110209}, it follows
that
\begin{equation}\label{eq-20110522}
   (a^q_1, \ldots, a^q_n) \succeq^\mathrm{ w} (b^q_1, \ldots, b^q_n)\quad\Longrightarrow\quad
       \sum^n_{i=1} a_i X_i \ge_\mathrm{ st} \sum^n_{i=1} b_i X_i
\end{equation}
when $p<0$, $q<0$ and $X_i^p$ has a log-concave density. Applying Lemma \ref{le-20110722}(iii)
with $g(x)=x^{\beta/q}$, we have
\[
   (a^q_1, \ldots, a^q_n) \succeq^\mathrm{ w} (b^q_1, \ldots, b^q_n)\quad\Longrightarrow\quad
       (a^\beta_1, \ldots, a^\beta_n) \succeq_\mathrm{ w} (b^\beta_1, \ldots, b^\beta_n)
\]
for any $q<0$ and $\beta>0$. Thus, \eqref{eq-20110522} can be deduced from Corollary
\ref{co-20110212-1}.

\item $\phi(x)=x^{1/q}$ ($q<0$) and $\psi(x)=\mathrm{e^\mathit{x}}$:\\
Conditions \eqref{eq-20110210-1} and \eqref{eq-20110210-2} are satisfied. From Theorem
\ref{th-20110209}, it follows that
\begin{equation}\label{eq-20110213-2}
   (a^q_1, \ldots, a^q_n) \succeq^\mathrm{ w} (b^q_1, \ldots, b^q_n)\quad\Longrightarrow\quad
       \sum^n_{i=1} a_i X_i \ge_\mathrm{ st} \sum^n_{i=1} b_i X_i
\end{equation}
when $q<0$ and $\log X_i$ has a log-concave density. It should be pointed out that
\eqref{eq-20110213-2} is implied by \eqref{eq-20110213-1} because applying Lemma
\ref{le-20110722}(iii) with $g(x)=q^{-1}\log x$ yields
\[
   (a^q_1, \ldots, a^q_n) \succeq^\mathrm{ w} (b^q_1, \ldots, b^q_n)\quad\Longrightarrow\quad
       (\log a_1, \ldots, \log a_n) \succeq_\mathrm{ w} (\log b_1, \ldots, \log b_n).
\]

\item $\phi(x)=\log (x+\mathrm{e})$ and $\psi(x)=x^{1/p}$ ($p\ge 2$):\\
Conditions \eqref{eq-20110210-2} and \eqref{eq-20110210-3} are satisfied.  Theorem
\ref{th-20110210} reduces to Corollary \ref{co-20110212-2} below, which can be deduced from Corollary \ref{th-Yu} by
observing $g(x)=(\log (x+\mathrm{e}))^q$ is increasing and concave on $\Re_+$ and applying Lemma \ref{le-20110722}(ii).
\end{itemize}
\end{remark}

\begin{corollary}\label{co-20110212-1}
Let $p<0$ and $0<q<1$ with $p^{-1}+q^{-1} = 1$, and let $X_1,X_2, \ldots, X_n$ be i.i.d. random
variables with density function $f$ on $\Re_+$. If $X^p_1$ has a log-concave density function,
then, for $\mathbf{a}, \mathbf{b}\in \Re_+^n$,
\[
   (a^q_1, \ldots, a^q_n) \succeq_\mathrm{ w} (b^q_1, \ldots, b^q_n)\quad\Longrightarrow\quad
       \sum^n_{i=1} a_i X_i \ge_\mathrm{ st} \sum^n_{i=1} b_i X_i.
\]
\end{corollary}

\begin{corollary}[(\cite{8})]\label{th-Karlin} Let $0<p<1$ and $q<0$ with $p^{-1}+q^{-1}=1$, and let $X_1,X_2,
\ldots, X_n$ be i.i.d. random variables with density function $f$ on $\Re_+$. If $X^p_1$ has a
log-concave density function, then
\[
   (a^q_1, \ldots, a^q_n) \succeq^\mathrm{ w} (b^q_1, \ldots, b^q_n)\quad\Longrightarrow\quad
       \sum^n_{i=1} a_i X_i \ge_\mathrm{ st} \sum^n_{i=1} b_i X_i.
\]
\end{corollary}

\begin{corollary}[(\cite{22})]\label{th-Yu} Let $p>1$ and $q>1$ with $p^{-1}+q^{-1}=1$, and let $X_1,X_2, \ldots, X_n$ be i.i.d. random variables with
density function $f$ on $\Re_+$. If $X^p_1$ has a log-concave density, then, for $\mathbf{a}, \mathbf{b}\in
\Re_+^n$,
\[
   (a^q_1, \ldots, a^q_n) \succeq^\mathrm{ w} (b^q_1, \ldots, b^q_n)\quad\Longrightarrow\quad
       \sum^n_{i=1} a_i X_i \le_\mathrm{ st} \sum^n_{i=1} b_i X_i.
\]
\end{corollary}

\begin{corollary}\label{co-20110212-2}
Let $p\ge 2$, and let $X_1,X_2, \ldots, X_n$ be i.i.d. random variables with density function $f$ on
$\Re_+$. If $X^p_1$ has a log-concave density function, then, for $\mathbf{a}, \mathbf{b}\in [1,\oo)^n$,
\[
   (\mathrm{e}^{a_1}, \ldots, \mathrm{e}^{a_n}) \succeq^\mathrm{ w} (\mathrm{e}^{b_1}, \ldots, \mathrm{e}^{b_n})\quad\Longrightarrow\quad
       \sum^n_{i=1} a_i X_i \le_\mathrm{ st} \sum^n_{i=1} b_i X_i.
\]
\end{corollary}

\begin{remark}\label{re-20110721}
In Theorems \ref{th-20110209} and \ref{th-20110210}, the functions $\phi$ and $\psi$ play an
independent role. From Remark \ref{re-20110213}, it is seen that, for $\psi(x)=\mathrm{e}^x$, $\phi(x)=\mathrm{e}^x$
is better than $\phi(x)=x^{1/q}$ ($q<0$). Two anonymous referees pointed out whether there is any
meaning of considering the best possible $\phi$ for a given $\psi$. This interesting question is
worth further investigation. We give a partial answer to this question.

For given two pairs $(\phi_1, \psi)$ and $(\phi_2, \psi)$ satisfying the conditions of
Theorem \ref{th-20110209} (resp. Theorem \ref{th-20110210}), define $g(x)=\phi_2^{-1}\circ \phi_1(x)$
(resp. $g(x)=-\phi_2^{-1}\circ \phi_1(x)$). By Lemma \ref{le-20110722},
$\phi_2$ is better than $\phi_1$ if either one of the following conditions holds:

\begin{enumerate}[(iii)]
  \item[(i)] $\phi_1$ and $\phi_2$ are increasing, and $g(x)$ is convex;

  \item[(ii)] $\phi_1$ is increasing, $\phi_2$ is decreasing, and $g(x)$ is concave;

  \item[(iii)] $\phi_1$ is decreasing, $\phi_2$ is increasing, and $g(x)$ is convex;

  \item[(iv)] $\phi_1$ and $\phi_2$ are decreasing, and $g(x)$ is concave.
\end{enumerate}
For example, denote $\psi_p(x)=x^{1/p}$ and $\phi_q(x)=x^{1/q}$ with $(p, q)\in A_1$ (resp. $A_2$,
$A_3$), where the $A_i$'s are defined in Remark \ref{re-20110213}. Fix $\psi_p(x)$, and choose
$q_\ast\in \Re$ such that $p^{-1}+q_{\ast}^{-1}=1$ and, hence, $(p, q_\ast)\in A_1$ (resp. $A_2$,
$A_3$). Define
\[
    g(x)=\phi_{q_\ast}^{-1}\circ \phi_q(x)=x^{q_\ast/q},\qquad x\in\Re_+.
\]
It is easy to see that, for $(p, q)\in A_1$ (resp. $A_2$, $A_3$), $g(x)$ is convex
(resp. concave, concave). Thus, for fixed $\psi_p$ with $(p,q)\in A_i$ ($i=1,2,3$),
$\phi_{q_\ast}$ is better than $\phi_q$.
\end{remark}

\section{Non-i.i.d. nonnegative random variables}\label{s4}

Before we prove the main results of this section, we give four lemmas. In the proofs of Theorems
\ref{th-20110209} and \ref{th-20110210}, we use an important fact that a permutation symmetric and
log-concave function is Schur-concave. In Lemma \ref{le-2011-1-19} below, a sufficient condition
is given for a log-concave function on $\Re_+^2$ to be Schur-concave on
$\mathscr{D}_+^2=\{(x_1,x_2)\dvtx x_1\le x_2, (x_1,x_2)\in\Re_+^2\}$. Lemma \ref{le-2011-1-19} plays a
key role in the proofs of Lemmas \ref{le-20110211-1} and \ref{le-20110211-2}.

\begin{lemma}\label{le-2011-1-19}
If $h(x_1, x_2)$ is log-concave on $\Re_+^2$ and
\[
   h\bigl(x_{(2)}, x_{(1)}\bigr) \ge h\bigl(x_{(1)}, x_{(2)}\bigr)\qquad \mbox{for all } (x_1,x_2)\in \Re_+^2,
\]
then
\[
   (x_1,x_2)\preceq_\mathrm{ m} (y_1, y_2)\quad\Longrightarrow\quad h\bigl(x_{(1)},x_{(2)}\bigr) \ge h\bigl(y_{(1)},y_{(2)}\bigr).
\]
\end{lemma}

\begin{pf} Suppose that $(x_1,x_2)\preceq_\mathrm{ m} (y_1, y_2)$. Then there exists $\alpha\in [1/2,1]$
such that
\[
   x_{(1)}= \alpha y_{(1)} + {\overline \alpha} y_{(2)},\qquad x_{(2)} =\alpha y_{(2)} + {\overline\alpha} y_{(1)}
\]
with ${\overline\alpha}=1-\alpha$. So,
\begin{eqnarray*}
  \log h\bigl(x_{(1)},x_{(2)}\bigr) &=&  \log h\bigl(\alpha y_{(1)} + {\overline \alpha} y_{(2)},\alpha y_{(2)}
                 + {\overline\alpha} y_{(1)}\bigr)\\
    &=& \log h\bigl(\alpha \bigl(y_{(1)}, y_{(2)}\bigr) + {\overline\alpha} \bigl(y_{(2)}, y_{(1)}\bigr)\bigr)\\
    & \ge & \alpha \log h\bigl(y_{(1)},y_{(2)}\bigr) + {\overline \alpha}  \log h\bigl(y_{(2)},y_{(1)}\bigr)\\
    & \ge &   \log h\bigl(y_{(1)},y_{(2)}\bigr),
\end{eqnarray*}
where the first inequality follows from the log-concavity of $h$.
\end{pf}

\begin{lemma}[(\cite{19})]\label{le-20110718}  Let $X$ and $Y$ be two independent random variables. Then
$X\ge_\mathrm{ lr} Y$ if and only if $g(X,Y)\ge_\mathrm{ st} g(Y,X)$ for all $g\in \mathscr{C}_\mathrm{ lr}$,
where
\[
   \mathscr{C}_\mathrm{ lr} =\{ g(x,y)\dvtx g(x, y)\ge g(y,x),\forall x\ge y\}.
\]
\end{lemma}

\begin{lemma}\label{le-20110211-1}
Let $X_1$ and $X_2$ be independent nonnegative random variables satisfying
\[
   X_1\ge_\mathrm{ lr} X_2,
\]
and let $\phi, \psi\dvtx \Re_+\to \Re_+$ be two twice continuously differentiable and strictly
monotone  functions such that \eqref{eq-20110210-1} and \eqref{eq-20110210-2} hold. If
$\psi^{-1}(X_i)$ has a log-concave density for each $i$, then
\[
    (\phi^{-1}(a_1), \phi^{-1}(a_2) )\succeq_\mathrm{ m}  (\phi^{-1}(b_1),
      \phi^{-1}(b_2) )\quad \Longrightarrow\quad a_{(2)} X_1+a_{(1)} X_2 \ge_\mathrm{ st} b_{(2)} X_1 + b_{(1)} X_2.
\]
\end{lemma}

\begin{pf} From the proof of Theorem \ref{th-20110209}, it follows that, for fixed $t\ge 0$,
\[
   h(c_1, c_2) = \P\bigl(\phi(c_1) X_1 + \phi(c_2) X_2 \le t \bigr),\qquad
      \eta(c_1, c_2) = \P\bigl(\phi(c_1) X_2 + \phi(c_2) X_1 \le t \bigr)
\]
are log-concave in $(c_1, c_2)\in \Re_+^2$ under conditions \eqref{eq-20110210-1} and
\eqref{eq-20110210-2}.

(1) Suppose that $\phi$ is decreasing. By Lemma \ref{le-20110718}, it follows that
\[
    \phi\bigl(c_{(2)}\bigr) X_1 + \phi\bigl(c_{(1)}\bigr) X_2 \le_\mathrm{ st}   \phi\bigl(c_{(1)}\bigr) X_1 + \phi\bigl(c_{(2)}\bigr) X_2,
\]
that is,
\[
    h\bigl(c_{(2)}, c_{(1)}\bigr) \ge h \bigl(c_{(1)}, c_{(2)}\bigr),\qquad (c_1,c_2)\in\Re_+^2.
\]
Then, by Lemma \ref{le-2011-1-19},
\[
   \phi\bigl(c_{(1)}\bigr) X_1 + \phi \bigl(c_{(2)}\bigr) X_2 \le_\mathrm{ st} \phi\bigl(d_{(1)}\bigr) X_1 + \phi \bigl(d_{(2)}\bigr) X_2
\]
whenever $(c_1,c_2), (d_1,d_2)\in\Re_+^2$ such that $(c_1, c_2)\preceq_\mathrm{ m} (d_1,d_2)$. Setting
$(c_1,c_2)=(\phi^{-1}(b_1)$, $\phi^{-1}(b_2))$ and $(d_1,d_2)=(\phi^{-1}(a_1), \phi^{-1}(a_2))$, it
follows that $ b_{(2)} X_1 + b_{(1)} X_2\le_\mathrm{ st} a_{(2)} X_1+a_{(1)} X_2$ since $\phi$ is deceasing.

(2) Suppose that $\phi$ is increasing. Again, by Lemma \ref{le-20110718}, it follows that
\[
    \eta\bigl(c_{(2)}, c_{(1)}\bigr) \ge \eta \bigl(c_{(1)}, c_{(2)}\bigr),\qquad (c_1,c_2)\in\Re_+^2.
\]
Then, by Lemma \ref{le-2011-1-19},
\[
   \phi\bigl(c_{(2)}\bigr) X_1 + \phi \bigl(c_{(1)}\bigr) X_2 \le_\mathrm{ st} \phi\bigl(d_{(2)}\bigr) X_1 + \phi \bigl(d_{(1)}\bigr) X_2
\]
whenever $(c_1,c_2), (d_1,d_2)\in\Re_+^2$ such that $(c_1, c_2)\preceq_\mathrm{ m} (d_1,d_2)$. This implies
$ b_{(2)} X_1 + b_{(1)} X_2\le_\mathrm{ st} a_{(2)} X_1+a_{(1)} X_2$ since $\phi$ is increasing.
This completes the proof of the lemma.
\end{pf}

\begin{lemma}\label{le-20110211-2}
Let $X_1$ and $X_2$ be independent nonnegative random variables satisfying
\[
   X_1\ge_\mathrm{ lr} X_2,
\]
and let $\phi, \psi\dvtx \Re_+\to \Re_+$ be two twice continuously differentiable and strictly
monotone functions such that \eqref{eq-20110210-2} and \eqref{eq-20110210-3} hold. If
$\psi^{-1}(X_i)$ has a log-concave density for each $i$, then
\[
    (\phi^{-1}(a_1), \phi^{-1}(a_2) )\succeq_\mathrm{ m}  (\phi^{-1}(b_1),
      \phi^{-1}(b_2) ) \quad\Longrightarrow \quad a_{(1)} X_1+a_{(2)} X_2 \le_\mathrm{ st} b_{(1)} X_1 + b_{(2)} X_2.
\]
\end{lemma}

\begin{pf} From the proof of Theorem \ref{th-20110210}, it follows that, for fixed $t\ge 0$,
\[
   h(c_1, c_2) = \P\bigl(\phi(c_1) X_1 + \phi(c_2) X_2 > t \bigr),\qquad
     \eta (c_1, c_2) = \P\bigl(\phi(c_2) X_1 + \phi(c_1) X_2 > t \bigr)
\]
are log-concave in $(c_1, c_2)\in \Re_+^2$ under conditions \eqref{eq-20110210-2} and
\eqref{eq-20110210-3}. The rest of the proof is similar to that of Lemma \ref{le-20110211-1} and
is, hence, omitted.
\end{pf}

Now, we are ready to present the following two main results, Theorems \ref{th-20110211-1} and
\ref{th-20110211-2}.

\begin{theorem}\label{th-20110211-1}
Let $X_1, X_2, \ldots, X_n$ be independent nonnegative random variables satisfying
\[
   X_1\ge_\mathrm{ lr} X_2\ge_\mathrm{ lr} \cdots \ge_\mathrm{ lr} X_n,
\]
and let $\phi, \psi\dvtx \Re_+\to \Re_+$ be two twice continuously differentiable and strictly
monotone functions such that \eqref{eq-20110210-1} and \eqref{eq-20110210-2} hold. Assume that
$\psi^{-1}(X_i)$ has a log-concave density for each $i$. If
\begin{equation}\label{eq-20110211-1}
      (\phi^{-1}(a_1),\ldots, \phi^{-1}(a_n) )\succeq_\mathrm{ m}  (\phi^{-1}(b_1), \ldots,
      \phi^{-1}(b_n) ),
\end{equation}
then
\begin{equation}\label{eq-20110720-1}
   \sum^n_{i=1} a_{(n-i+1)} X_i \ge_\mathrm{st} \sum^n_{i=1} b_{(n-i+1)} X_i.
\end{equation}
Moreover, if $\phi$ is increasing (decreasing), the majorization order in \eqref{eq-20110211-1}
can be replaced by the submajorization (supmajorization) order.
\end{theorem}

\begin{pf} By the nature of the supmajorization and submajorization orders (see the proof of
Theorem \ref{th-20110209}), it suffices to prove that \eqref{eq-20110211-1} implies
\eqref{eq-20110720-1}. Suppose that \eqref{eq-20110211-1} holds. Then, by Lem\-ma~2.B.1 in \cite{14},
 there exists a finite number of vectors $\phi^{-1}(\mathbf{c}^j_{()}):=(\phi^{-1}(c^j_{(1)}),\ldots, \phi^{-1}(c^j_{(n)}))\in \Re_+^n$, $j=1,\ldots, N$, such
that
\begin{eqnarray*}
   \bigl(\phi^{-1}\bigl(a_{(1)}\bigr),\ldots, \phi^{-1}\bigl(a_{(n)}\bigr) \bigr) &=& \phi^{-1}\bigl(\mathbf{c}^1_{()}\bigr) \preceq_\mathrm{m}
        \phi^{-1}\bigl(\mathbf{c}^2_{()}\bigr) \preceq_\mathrm{m} \cdots  \preceq_\mathrm{m} \phi^{-1}\bigl(\mathbf{c}^N_{()}\bigr)\\
      &=&   \bigl(\phi^{-1}\bigl(b_{(1)}\bigr), \ldots, \phi^{-1}\bigl(b_{(n)}\bigr) \bigr),
\end{eqnarray*}
where $\mathbf{c}^k_{()}=(c^k_{(1)}, \ldots, c^k_{(n)})$, the ordered vector of $\mathbf{c}^k=(c^k_1,
\ldots, c^k_n)\in \Re_+^n$, and $\mathbf{c}^k_{()}$ and $\mathbf{c}^{k+1}_{()}$ differ only in two
coordinates for each $k$. Therefore, the desired result now follows from Lemma \ref{le-20110211-1}
and the fact that the usual stochastic order is closed under convolution. This completes the
proof.
\end{pf}

Similarly, we can prove the next result by using Lemma \ref{le-20110211-2}.

\begin{theorem}\label{th-20110211-2}
Let $X_1, X_2, \ldots, X_n$ be independent nonnegative random variables satisfying
\[
    X_1\ge_\mathrm{lr} X_2\ge_\mathrm{lr} \cdots\ge_\mathrm{lr} X_n,
\]
and let $\phi, \psi\dvtx \Re_+\to \Re_+$ be two twice continuously differentiable and strictly
monotone functions such that \eqref{eq-20110210-2} and \eqref{eq-20110210-3} hold. Assume that
$\psi^{-1}(X_i)$ has a log-concave density for each $i$. If
\begin{equation}\label{eq-20110211-2}
    (\phi^{-1}(a_1), \ldots,\phi^{-1}(a_n) )\succeq_\mathrm{m}  (\phi^{-1}(b_1),\ldots,
      \phi^{-1}(b_n) ),
\end{equation}
then
\[
     \sum^n_{i=1} a_{(i)} X_i \le_\mathrm{st} \sum^n_{i=1} b_{(i)} X_i.
\]
Moreover, if $\phi$ is increasing (decreasing), the majorization order in \eqref{eq-20110211-2}
can be replaced by the supmajorization (submajorization) order.
\end{theorem}

Combining Theorems \ref{th-20110211-1} and \ref{th-20110211-2} and Remark \ref{re-20110213}, we
have the following corollaries, which extend some results in Section \ref{s3} from i.i.d. to non-i.i.d.
nonnegative random variables.

\begin{corollary}[(\cite{21})]\label{co-20110213-1}
 Let $X_1,X_2,\ldots,X_n$ be independent nonnegative random variables
satisfying $X_1\ge_\mathrm{lr} X_2\ge_\mathrm{lr} \cdots \ge_\mathrm{lr} X_n$. If $\log X_i$ has a log-concave
density for each $i$, then,
\[
   (\log a_1,\ldots, \log a_n)\succeq_\mathrm{w} (\log b_1,\ldots, \log b_n)
      \quad\Longrightarrow \quad\sum_{i=1}^n a_{(n-i+1)} X_i \ge_\mathrm{st} \sum_{i=1}^n b_{(n-i+1)} X_i.
\]
\end{corollary}

\begin{corollary}\label{co-20110118}
Let $p>1$ and $q>1$ with $p^{-1}+q^{-1}=1$, and let $X_1,X_2, \ldots, X_n$ be independent
nonnegative random variables such that $X_1\ge_\mathrm{lr} X_2\ge_\mathrm{lr} \cdots \ge_\mathrm{lr} X_n$.
If $X^p_i$ has a log-concave density function for each $i$, then, for $\mathbf{a}, \mathbf{b}\in
\Re_+^n$,
\[
    (a^q_1, \ldots, a^q_n) \preceq^\mathrm{w} (b^q_1, \ldots, b^q_n)\quad \Longrightarrow\quad
        \sum^n_{i=1} a_{(i)} X_i \ge_\mathrm{st} \sum^n_{i=1} b_{(i)} X_i.
\]
\end{corollary}

\begin{corollary}\label{co-20110119}
Let $p\in (0,1)$ and $q<0$ with $p^{-1}+q^{-1}=1$, and let $X_1,X_2, \ldots, X_n$ be  independent
nonnegative random variables such that $X_1\ge_\mathrm{lr} X_2\ge_\mathrm{lr} \cdots \ge_\mathrm{lr} X_n$.
If $X^p_i$ has a log-concave density function for each $i$, then, for $\mathbf{a}, \mathbf{b}\in
\Re_+^n$,
\[
   (a^q_1, \ldots, a^q_n) \succeq^\mathrm{w} (b^q_1, \ldots, b^q_n)\quad\Longrightarrow\quad
       \sum^n_{i=1} a_{(n-i+1)} X_i \ge_\mathrm{st} \sum^n_{i=1} b_{(n-i+1)} X_i.
\]
\end{corollary}

\begin{corollary}\label{co-20110213-2}
Let $p<0$ and $0<q<1$ with $p^{-1}+q^{-1} = 1$, and let $X_1,X_2, \ldots, X_n$ be independent
nonnegative random variables such that $X_1\ge_\mathrm{lr} X_2\ge_\mathrm{lr} \cdots \ge_\mathrm{lr} X_n$.
If $X^p_i$ has a log-concave density function for each $i$, then, for $\mathbf{a}, \mathbf{b}\in
\Re_+^n$,
\[
   (a^q_1, \ldots, a^q_n) \succeq_\mathrm{w} (b^q_1, \ldots, b^q_n)\quad\Longrightarrow\quad
       \sum^n_{i=1} a_{(n-i+1)} X_i \ge_\mathrm{st} \sum^n_{i=1} b_{(n-i+1)} X_i.
\]
\end{corollary}

Finally, we give an example to which Corollaries \ref{co-20110213-1}--\ref{co-20110213-2} can be applied.

\begin{example}
Let $X$ be a nonnegative random variable having the generalized gamma distribution $F_{p,\alpha,
\lambda}$ with density function
\[
   f_{p,\alpha, \lambda}(x) =\frac {p\lambda^\alpha}{\Gamma(\alpha)} x^{\alpha p-1} \mathrm{e}^{-\lambda x^p}, \qquad x>0,
\]
where $p>0$, $\alpha>0$ and $\lambda>0$ (see \cite{7,16}). This class of
distributions includes the Weibull $(\alpha=1)$, gamma $(p=1)$ and the generalized Rayleigh $(p=2)$
distributions as special cases. It is easy to see that
\begin{itemize}
  \item for $\alpha\ge 1$, $X^p$ has a log-concave density;
  \item for $0<\alpha<1$, $X^{\alpha p}$ has a log-concave density;
  \item for $0<\alpha_1 \le\alpha_2$, $F_{p,\alpha_1,\lambda}\le_\mathrm{lr} F_{p,\alpha_2,\lambda}$;
  \item for $0<\lambda_1\le \lambda_2$, $F_{p,\alpha,\lambda_2}\le_\mathrm{lr} F_{p,\alpha,\lambda_1}$;
  \item $\log X$ has a log-concave density.
\end{itemize}

Corollaries \ref{co-20110118} and \ref{co-20110119} can be applied to the above generalized gamma
distribution. For example, let $X_1, X_2, \ldots, X_n$ be independent nonnegative random variables
having distributions $F_{p,\alpha, \lambda_1}, F_{p,\alpha, \lambda_2}, \ldots, F_{p,\alpha,
\lambda_n}$ with $0<\lambda_1\le \lambda_2\le \cdots\le \lambda_n$ and $\alpha\ge 1$, or having
distributions $F_{p,\alpha_1, \lambda}, F_{p,\alpha_2, \lambda}, \ldots, F_{p,\alpha_n, \lambda}$
with $\alpha_1\ge \alpha_2\ge \cdots\ge \alpha_n\ge 1$. Then Corollaries \ref{co-20110118} and
\ref{co-20110119} hold for $p>1$ and $p<1$, respectively.
\end{example}

\section*{Acknowledgements}

We thank the referees for comments on a previous draft of the paper. The comments led us to
significantly improve the paper. T. Hu supported by the NNSF of China (Nos. 11071232, 70821001),  and the National Basic Research
 Program of China (973 Program, Grant No.~2007CB814901).


\printhistory

\end{document}